\documentclass[12pt,reqno]{amsart}
\usepackage[margin=0.95in]{geometry}
\usepackage{amsmath,amssymb,amsthm,mathtools}
\usepackage{enumitem}
\usepackage{microtype}
\usepackage{mathrsfs}
\usepackage[hidelinks]{hyperref}

\numberwithin{equation}{section}

\newcommand{\Irr}{\operatorname{Irr}}
\newcommand{\cod}{\operatorname{cod}}
\newcommand{\cd}{\operatorname{cd}}

\newcommand{\Sol}{\operatorname{Sol}}
\newcommand{\Aut}{\operatorname{Aut}}
\newcommand{\GL}{\operatorname{GL}}
\newcommand{\SL}{\operatorname{SL}}
\newcommand{\Sp}{\operatorname{Sp}}

\newcommand{\F}{\operatorname{F}}
\newcommand{\PhiG}{\operatorname{\Phi}}

\newcommand{\OO}{\operatorname{O}}
\newcommand{\dl}{\operatorname{dl}}

\newcommand{\C}{\operatorname{C}}

\newtheorem*{thmA}{Lemma A}
\newtheorem*{thmB}{Theorem B}
\newtheorem*{thmC}{Theorem C}

\newtheorem{theorem}{Theorem}[section]
\newtheorem{lemma}[theorem]{Lemma}

\newtheorem{proposition}[theorem]{Proposition}
\theoremstyle{remark}

\title{Solvable groups with square-free character degree--codegree gcds}
\author{KARAM ALDAHLEH}
\address{Mathematical Institute, University of Oxford, Oxford, OX2 6GG, UK.}
\email{karam.aldahleh@hertford.ox.ac.uk}
\author{YONG YANG}
\address{Department of Mathematics, Texas State University, San Marcos, TX 78666, USA.}
\email{yang@txstate.edu}


\begin{document}

\begin{abstract}
For an irreducible character $\chi$ of a finite group $G$, its codegree is
$\cod(\chi)=|G:\ker\chi|/\chi(1)$. We study finite solvable groups $G$
for which $\gcd(\chi(1),\cod(\chi))$ is square-free for every
irreducible character $\chi$. We characterize the nilpotent groups with this
property and prove that the Fitting height and derived length are at
most $5$ and $6$, respectively. Both bounds are sharp. Writing $\F(G)$
for the Fitting subgroup, we also show that the quotients $G/\F(G)$
have unbounded Pr\"ufer rank.
\end{abstract}

\maketitle

\section{Introduction}\label{sec:introduction}

All groups are finite, and all characters are complex characters.
We write $\Irr(G)$ for the set of irreducible characters of $G$ and
$\cd(G)=\{\chi(1):\chi\in\Irr(G)\}$ for its set of character degrees.
For $\chi\in\Irr(G)$, the codegree is
\[
 \cod(\chi)=\frac{|G:\ker\chi|}{\chi(1)}.
\]
This invariant was introduced in \cite{CodegreeDefn}.
We denote the Fitting and Frattini subgroups by $\F(G)$ and
$\PhiG(G)$, respectively, and the largest normal $p$-subgroup by
$\OO_p(G)$. For a solvable group $G$, its Fitting height and derived
length are denoted by $h(G)$ and $\dl(G)$. We use
$G^{(0)}=G$ and $G^{(i+1)}=[G^{(i)},G^{(i)}]$ for the derived series.
$Z(G)$ denotes the center of $G$. We write $C_m$ for the cyclic
group of order $m$ and $S_n$ for the symmetric group of degree $n$.
The minimum number of generators of a group $H$ is $d(H)$; its
Pr\"ufer rank is $\max\{d(K):K\le H\}$. Throughout, rank means
Pr\"ufer rank.

The common divisors of $\chi(1)$ and $\cod(\chi)$ impose strong
restrictions on $G$. Groups for which these two integers are coprime
for every irreducible character are called $\mathscr H$-groups and
were classified in \cite{Liang}. More generally, for a prime $p$, an
$\mathscr H_p$-group satisfies
$p\nmid\gcd(\chi(1),\cod(\chi))$ for every $\chi\in\Irr(G)$;
these groups were studied in \cite{Qian-Zeng}. Related conditions
appear in \cite{Chen-Lv,Gao-Wang-Chen,Lv-Yang-Qian}.

In \cite[Problem~5.7]{Qian}, the following condition was proposed:
\begin{equation}\label{eq:sf}
 \gcd(\chi(1),\cod(\chi))\text{ is square-free for every }
 \chi\in\Irr(G).
\end{equation}
The problem asks for bounds on the Fitting height and the rank of
$G/\F(G)$ when $G$ is solvable, and for a description of the
nonsolvable groups satisfying \eqref{eq:sf}. Let $\Sol(G)$ denote
the largest solvable normal subgroup of $G$. The possible quotients
$G/\Sol(G)$ in the nonsolvable case were determined in
\cite[Theorem~A]{AldahlehKapplerMakurYang}. Here we consider the
solvable case.

For a nilpotent group, every irreducible character degree divides
its codegree \cite{Gagola-Lewis}. Thus \eqref{eq:sf} is equivalent
to all character degrees being square-free. This gives the following description.

\begin{thmA}\label{lem:nilpotent}
A nilpotent group $G$ satisfies \eqref{eq:sf} if and only if, for every
prime $p$ dividing $|G|$, its Sylow $p$-subgroup $P$ satisfies at least
one of the following conditions:
\begin{enumerate}[label=\textup{(\roman*)}]
\item $P$ is abelian;
\item $P$ has an abelian normal subgroup of index $p$;
\item $|P:Z(P)|=p^3$, and $P$ has an abelian normal subgroup of
index at most $p^2$.
\end{enumerate}
\end{thmA}

Solvable $\mathscr H$-groups have Fitting height at most $3$
\cite[Theorem~A]{Liang}. Under \eqref{eq:sf}, the best possible bound
is $5$.

\begin{thmB}\label{thm:fitting}
Let $G$ be a solvable group satisfying \eqref{eq:sf}. Then $h(G)\le5$.
\end{thmB}

Lemma A implies that the nilpotent factors in the Fitting series are
metabelian. Together with Theorem B, this gives $\dl(G)\le10$.
The following result improves this bound to $6$.

\begin{thmC}\label{thm:derived}
Let $G$ be a solvable group satisfying \eqref{eq:sf}. Then $\dl(G)\le6$.
\end{thmC}

Both bounds are attained by the groups constructed in
Section~\ref{sec:examples}. In contrast, the rank of $G/\F(G)$ is
unbounded: the direct product $S_3^n$ satisfies \eqref{eq:sf}, and
$S_3^n/\F(S_3^n)\cong C_2^n$. 

Section~\ref{sec:fitting} proves Lemma A and Theorem B.
Section~\ref{sec:derived} first establishes the required bounds for
completely reducible linear groups of degree at most $3$. These bounds
are then applied to primitive affine sections and finally to the
irreducible constituents of $\F(G)/\PhiG(G)$.

We use $\chi_H$ for restriction, $\theta^G$ for induction, and
$\langle\ ,\ \rangle_H$ for the usual character inner product on $H$.
If $N\lhd G$ and $\theta\in\Irr(N)$, then $\Irr(G\mid\theta)$
denotes the set of irreducible characters of $G$ whose restrictions
to $N$ contain $\theta$. For a group $H$ acting on a character $\theta$,
$I_H(\theta)$ denotes its stabilizer. For a positive integer $m$,
$v_p(m)$ is the exponent of $p$ in $m$, and $m_p=p^{v_p(m)}$.
The dual and exterior powers of a vector space $V$ over its base
field are denoted by $V^*$ and $\bigwedge^jV$.

\section{Nilpotent groups and Fitting height}\label{sec:fitting}

\begin{lemma}\label{lem:inheritance}
If $G$ satisfies \eqref{eq:sf}, then every quotient of every subnormal
subgroup of $G$ satisfies \eqref{eq:sf}.
\end{lemma}

\begin{proof}
Inflation preserves character degree and codegree, so the condition
passes to quotients. Let $H$ be subnormal in $G$ and let
$\theta\in\Irr(H)$. Choose $\chi\in\Irr(G)$ whose restriction to $H$
contains $\theta$. Clifford theory along a subnormal series gives
$\theta(1)\mid\chi(1)$, and \cite[Lemma~2.1(2)]{Liang} gives
$\cod(\theta)\mid\cod(\chi)$. Hence
\[
 \gcd(\theta(1),\cod(\theta))
 \mid\gcd(\chi(1),\cod(\chi)),
\]
as required.
\end{proof}

\begin{proof}[Proof of Lemma A]
For nilpotent $G$, we have $\chi(1)\mid\cod(\chi)$ for every
$\chi\in\Irr(G)$. Since $G$ is the direct product of its Sylow
subgroups, \eqref{eq:sf} holds if and only if
$\cd(P)\subseteq\{1,p\}$ for each Sylow $p$-subgroup $P$.

Suppose that $P$ is nonabelian and $\cd(P)=\{1,p\}$. Choose
$N\lhd P$ maximal subject to $H=P/N$ being nonabelian. Every
nontrivial normal subgroup of $H$ contains $H'$. In particular,
a subgroup of order $p$ in $Z(H)$ contains $H'$, so $H'$ is
central of order $p$. A nonlinear character $\alpha\in\Irr(H)$
has degree $p$ and is faithful. Indeed, a nontrivial kernel would
contain $H'$.

Let $\rho$ afford $\alpha$. For $g\notin Z(H)$, choose $t\in H$
with $[g,t]\ne1$. Since $[g,t]\in H'$, its image under $\rho$ is
$\zeta I$ for some $\zeta\ne1$. Thus
\[
 \alpha(g)=\alpha(t^{-1}gt)=\zeta\alpha(g),
\]
so $\alpha(g)=0$. On the other hand, $|\alpha(g)|=p$ for $g\in Z(H)$. Character
orthogonality gives $|H|=p^2|Z(H)|$. Let $Q$ be the inverse image
of $Z(H)$ in $P$, and let $\beta$ be the inflation of $\alpha$ to
$P$. Then $|P:Q|=p^2$, the character $\beta$ vanishes outside $Q$,
and $\beta_Q=p\lambda$ for a linear character $\lambda$ of $Q$.

The subgroup $Q$ is abelian. Otherwise choose a nonlinear
$\theta\in\Irr(Q)$, together with
$\chi\in\Irr(P\mid\theta)$ and
$\eta\in\Irr(P\mid\theta\lambda)$. Both characters are nonlinear,
so $\chi(1)=\eta(1)=p$. Since $\theta(1)$ is a power of $p$
greater than $1$, we have $p\le\theta(1)\le\chi(1)=p$.
Clifford theory gives
$\chi_Q=\theta$ and $\eta_Q=\theta\lambda$.
Consequently,
\[
 \langle\chi\beta,\eta\rangle_P
 =\frac{p}{|P|}\sum_{g\in Q}(\theta\lambda)(g)
                  \overline{(\theta\lambda)(g)}
 =\frac{p|Q|}{|P|}=\frac1p,
\]
contrary to the integrality of character inner products.

Suppose now that $P$ has no abelian subgroup of index $p$.
Then $\C_P(Q)=Q$, since any element of $\C_P(Q)\setminus Q$
would generate with $Q$ an abelian subgroup of index at most $p$.
Let $\ell$ be the number of orbits of $Q$ on $P$ under conjugation.
Orbit counting and column orthogonality give
\[
 \ell=\frac1{|Q|}\sum_{g\in Q}|\C_P(g)|
 =\sum_{\psi\in\Irr(P)}\langle\psi_Q,\psi_Q\rangle_Q.
\]
Since $Q$ is abelian, the multiplicities of the linear constituents
of $\psi_Q$ give
$\langle\psi_Q,\psi_Q\rangle_Q\ge\psi(1)$. Therefore
\[
 \ell\ge\sum_{\psi\in\Irr(P)}\psi(1)
 \ge\frac1p\sum_{\psi\in\Irr(P)}\psi(1)^2=\frac{|P|}{p}.
\]
The singleton orbits are exactly the elements of $Q$. If every other
orbit had size at least $p^2$, then
\[
 \ell\le |Q|+\frac{|P|-|Q|}{p^2}
 =|Q|\left(2-\frac1{p^2}\right)<p|Q|=\frac{|P|}{p},
\]
a contradiction. Hence some $x\in P\setminus Q$ satisfies
$|Q:\C_Q(x)|=p$.

Let $A=\langle Q,x\rangle$ and $B=Z(A)$. Since $P/Q$ is abelian,
$A\lhd P$, and therefore $B\lhd P$. Also
$B\le\C_P(Q)=Q$, whence $B=\C_Q(x)$ and $|Q:B|=p$.
The group $A$ is nonabelian. Suppose that $[B,P]\ne1$, and choose
$1\ne u\in A'$ and $1\ne v\in[B,P]$. Both elements lie in $Q$.
There is a linear character $\mu$ of $Q$ with
$\mu(u)\ne1$ and $\mu(v)\ne1$: the two annihilators are proper
subgroups of $\Irr(Q)$ and cannot cover it.
Choose $\psi\in\Irr(P\mid\mu)$. Since $u\notin\ker\psi$,
the restriction $\psi_A$ has a nonlinear constituent. Since
$\psi(1)\le p$, this constituent has degree $p$ and equals $\psi_A$.
Every element of $B$
therefore acts as a scalar in a representation affording $\psi$.
This gives $[B,P]\le\ker\psi$, contrary to $\mu(v)\ne1$.
Thus $B\le Z(P)$. The reverse inclusion follows from
$Z(P)\le Q\le A$, so $B=Z(P)$ and
\[
 |P:Z(P)|=|P:Q||Q:B|=p^3.
\]
The abelian normal subgroup $Q$ has index $p^2$, as required.

Conversely, if $P$ is abelian or has an abelian normal subgroup of
index $p$, then $\cd(P)\subseteq\{1,p\}$; in the second case,
every character degree divides that index by
\cite[Theorem~6.15]{Isaacs/book}. If $|P:Z(P)|=p^3$, then
\[
 \chi(1)^2\le |P:Z(P)|=p^3
 \qquad(\chi\in\Irr(P)).
\]
As $\chi(1)$ is a power of $p$, it is either $1$ or $p$.
This proves the converse.
\end{proof}

Every nilpotent group satisfying \eqref{eq:sf} is metabelian by
Lemma A: each Sylow subgroup has an abelian normal subgroup with
abelian quotient. Lemma~\ref{lem:inheritance} therefore shows that
the nilpotent factors in the Fitting series of a solvable group
satisfying \eqref{eq:sf} are metabelian. In particular,
\[
 \dl(G)\le2h(G).
\]

For a group $H$ acting on a vector space $V$, we write $\C_H(v)$
for the stabilizer of $v$, $\C_H(V)$ for the kernel of the action,
and $\C_V(H)$ for the fixed subspace. A linear group is called
primitive if it is irreducible and preserves no direct-sum
decomposition into more than one nonzero subspace permuted by the
group.

\begin{lemma}\label{lem:affine-rank}
Let $V$ be a nonzero vector space over $\mathbb F_q$, where $q$ is
prime, and let $H$ act faithfully and irreducibly on $V$. Suppose
that $V\rtimes H$ is solvable and satisfies \eqref{eq:sf}. Then
$\OO_q(H)=1$. If $P$ is a nontrivial normal $p$-subgroup of $H$,
then $p\ne q$ and the following hold:
\begin{enumerate}[label=\textup{(\roman*)}]
\item For every $0\ne v\in V^*$,
\begin{equation}\label{eq:orbit}
 |\C_P(v)|\le p\qquad\text{or}\qquad|P:\C_P(v)|\le p.
\end{equation}
\item Every abelian subgroup $A$ of $P$ satisfies $d(A)\le3$.
If $d(A)=3$, then $A\cong C_p^3$.
\item We have $d(P)\le4$. If $P$ has an abelian subgroup of index
$p$, then $d(P)\le3$. Every irreducible $H$-composition factor of
$P/\PhiG(P)$ has dimension at most $3$.
\end{enumerate}
\end{lemma}

\begin{proof}
The fixed subspace of the normal $q$-subgroup $\OO_q(H)$ on $V$
is nonzero and $H$-invariant. Irreducibility and faithfulness
therefore give $\OO_q(H)=1$.

Let $G=V\rtimes H$. Identify $\Irr(V)$ with $V^*$ by fixing a
primitive $q$th root of unity. For $1\ne\lambda\in\Irr(V)$, let
$T=I_H(\lambda)$. Extend $\lambda$ to $VT$ by making $T$ act
trivially, and denote this extension by $\widetilde\lambda$.
Then $\chi=\widetilde\lambda^{\,G}$ is irreducible and
$\chi(1)=|H:T|$. The $H$-orbit of $\lambda$ spans $V^*$, so
$\ker\chi\cap V=1$. If $K=\ker\chi$, then $[K,V]=1$, whence
$K\le\C_G(V)=V$. Thus $\chi$ is faithful, and
\[
 \chi(1)=|H:T|,\qquad \cod(\chi)=|V||T|.
\]
For $p\ne q$, condition \eqref{eq:sf} gives
\[
 \min\{v_p(|H:T|),v_p(|T|)\}\le1.
\]
Since $P\lhd H$, we have
$|P:\C_P(\lambda)|=|PT:T|\mid|H:T|$ and
$\C_P(\lambda)\le T$. This proves \eqref{eq:orbit}.

Let $A\le P$ be abelian. The orbit condition also holds with $A$
in place of $P$, since stabilizer orders and orbit lengths for
$A$ divide the corresponding quantities for $P$. As $p\ne q$,
regard $V^*$ as an $A$-module and write
$V^*=U_1\oplus\cdots\oplus U_t$, where each $U_i$ is irreducible.
Let $K_i=\C_A(U_i)$. The group $A/K_i$ is cyclic
and acts freely on $U_i\setminus\{0\}$. Faithfulness gives
$\bigcap_iK_i=1$, and \eqref{eq:orbit} gives
\[
 |K_i|\le p\qquad\text{or}\qquad |A:K_i|\le p.
\]
If $|K_i|\le p$ for some $i$, then $d(A)\le2$, since $A/K_i$
is cyclic. Otherwise every proper $K_i$ has index $p$, so
$\PhiG(A)\le\bigcap_iK_i=1$ and $A\cong C_p^r$. If $r\ge4$,
choose distinct proper $K_i,K_j$ and nonzero $u_i\in U_i$,
$u_j\in U_j$. Then
\[
 \C_A(u_i+u_j)=K_i\cap K_j
\]
has index $p^2$ and order at least $p^2$, contrary to
\eqref{eq:orbit}. This proves \textup{(ii)}.

By Lemma~\ref{lem:inheritance}, the group $P$ satisfies
\eqref{eq:sf}. If $P$ is abelian, then $d(P)\le3$ by
\textup{(ii)}. Otherwise choose an abelian normal subgroup $A$
of index $p$ if one exists; if not, choose $A$ as in Lemma A\textup{(iii)}.
In either case $|P:A|\le p^2$. If $d(A)\le2$, then
\[
 d(P)\le d(A)+d(P/A)\le4,
\]
and this bound is at most $3$ when $|P:A|=p$.
If $d(A)=3$, then $A\cong C_p^3$. The subgroup $A$ is not central:
this follows from the nonabelianness of $P$ when $|P:A|=p$, and
from $|P:Z(P)|=p^3$ when $|P:A|=p^2$.
Thus $1\ne[A,P]\le A\cap\PhiG(P)$, and the image of $A$ in
$P/\PhiG(P)$ has dimension at most $2$. The quotient by this
image contributes at most one dimension if $|P:A|=p$, and at most
two if $|P:A|=p^2$. The stated generator bounds follow.

It remains to consider $d(P)=4$. The preceding bounds exclude
cases \textup{(i)} and \textup{(ii)} of Lemma A, so
$|P:Z(P)|=p^3$. The $H$-invariant subspace
\[
 Z(P)\PhiG(P)/\PhiG(P)
\]
of $P/\PhiG(P)$ is nonzero, since otherwise
$|P:Z(P)|\ge|P:\PhiG(P)|=p^4$. It is proper, since
$Z(P)\PhiG(P)=P$ would force $P$ to be generated by central
elements. Hence every composition factor has dimension at most $3$.
\end{proof}

\begin{proof}[Proof of Theorem B]
Let $G$ be a counterexample of minimal order. Quotients inherit
\eqref{eq:sf} by Lemma~\ref{lem:inheritance}, and factoring out
the Frattini subgroup preserves Fitting height
\cite[Chapter~III, Section~4]{DoerkHawkes}. Thus $\PhiG(G)=1$.
If $G$ had distinct minimal normal subgroups $N_1,N_2$, then it
would embed in $(G/N_1)\times(G/N_2)$, whose Fitting height is at
most $5$. Subgroups of a solvable group cannot have larger Fitting
height, so $G$ has a unique minimal normal subgroup $V$.
By \cite[Theorem~1.12]{ManzWolf}, we have $\F(G)=V$ and
$G=V\rtimes H$, where $V$ is elementary abelian and $H$ acts
faithfully and irreducibly on $V$. Hence
\begin{equation}\label{eq:height-split}
 h(G)=1+h(H).
\end{equation}

Let $F=\F(H)$. By \cite[Theorem~1.12]{ManzWolf}, the group
$F/\PhiG(H)$ is a faithful completely reducible $H/F$-module.
For each prime $p\mid|F|$, its $p$-component is
\[
 W_p=\frac{\OO_p(H)\PhiG(H)}{\PhiG(H)}
 \cong\frac{\OO_p(H)}{\OO_p(H)\cap\PhiG(H)}.
\]
Since $W_p$ is elementary abelian, it is a quotient of
$\OO_p(H)/\PhiG(\OO_p(H))$. Lemma~\ref{lem:affine-rank}
therefore shows that every irreducible constituent of $W_p$
has dimension at most $3$.

A solvable completely reducible linear group of degree at most $3$
has Fitting height at most $3$ by \cite[Theorem~C]{FrickNewman}.
The group $H/F$ embeds in the direct product of its images on the
irreducible constituents of $F/\PhiG(H)$. Each image has Fitting
height at most $3$, so $h(H/F)\le3$. Consequently $h(H)\le4$,
and \eqref{eq:height-split} gives $h(G)\le5$, a contradiction.
\end{proof}

\section{Derived length}\label{sec:derived}

We first bound the derived length of solvable linear groups of
degree at most $3$. These bounds will be applied to complements
in primitive solvable groups and then to the irreducible constituents
of $\F(G)/\PhiG(G)$.

\begin{lemma}\label{lem:linear-derived}
Let $r$ be a prime, and let $G\le\GL_d(\mathbb F_r)$ be solvable and
completely reducible. Suppose that $G$ satisfies \eqref{eq:sf}. Then
$\dl(G)\le3$ if $d\le2$, and $\dl(G)\le4$ if $d\le3$.
\end{lemma}

\begin{proof}
The images of $G$ on its irreducible constituents are quotients of $G$,
and $G$ embeds in their direct product. We may therefore assume that
$G$ is irreducible. The assertion is immediate when $d=1$. If $G$ is
imprimitive, then its blocks are one-dimensional. The kernel of the
action on these blocks is abelian, and the permutation image is a
subgroup of $S_d$. Thus $\dl(G)\le2$ when $d=2$, and
$\dl(G)\le3$ when $d=3$.

Suppose that $G$ is primitive, and choose a maximal solvable subgroup
$H$ of $\GL_d(\mathbb F_r)$ containing $G$. This overgroup is also
primitive. By \cite[Lemma~2.2]{Seress}, $H$ has normal subgroups
$A\le B\le C$ with the following properties. The group $A$ is cyclic,
and its linear span is a field $\mathbb F_{r^a}$ for some positive
integer $a$. Also $C=\C_H(A)$.
There is a positive integer $e$ with $d=ae$, the group $B$ has class
at most $2$, and $A\le Z(B)$. Moreover, $H/C$ is cyclic of order
dividing $a$, while $|B:A|=e^2$. When $e=1$, we have $B=C=A$.
When $e=d$, the group $A$ consists of scalars, $C=H$, and
\[
 B/A\cong C_d^2,\qquad H/B\le\Sp_2(d),
\]
where $\Sp_2(d)$ is the symplectic group on a $2$-dimensional
space over $\mathbb F_d$. The latter action is faithful and preserves
the nondegenerate
commutator form on $B/A$.

Since $d\in\{2,3\}$, either $e=1$ or $e=d$. In the first case
$H$ is metabelian. In the second case, set $N=G\cap B$. Then $N$
has class at most $2$ and
\[
 G/N\le\Sp_2(d)\cong
 \begin{cases}
  S_3,&d=2,\\
  \SL_2(3),&d=3.
 \end{cases}
\]
Consequently $\dl(G)\le d+2$. Suppose that equality holds.
Every proper subgroup of $S_3$ is abelian, and every proper subgroup of
$\SL_2(3)$ is metabelian. Thus $N$ is nonabelian and
$G/N=\Sp_2(d)$. Since $N\nleq A$, the subspace $NA/A$ of $B/A$
is nonzero and $G$-invariant. Since $G$ induces the full symplectic
group on $B/A$, this module is irreducible, and hence $NA/A=B/A$.
The commutator form on $B/A$ is nondegenerate, so
\[
 Z(N)=N\cap A,\qquad N/Z(N)\cong C_d^2.
\]
In particular, $Z(N)$ is cyclic. A faithful character of $Z(N)$ lies
under a unique faithful $\theta\in\Irr(N)$ of degree $d$. Since
$Z(N)$ is central in $G$, the character $\theta$ is $G$-invariant.
Both $S_3$ and $\SL_2(3)$ have trivial Schur multiplier, so $\theta$
extends to a character $\widehat\theta\in\Irr(G)$ by
\cite[Theorem~11.7]{Isaacs/book}. Each quotient has
an irreducible character $\tau$ of degree $d$. Inflating $\tau$ to
$G$, \cite[Corollary~6.17]{Isaacs/book} gives
\[
 \chi=\widehat\theta\tau\in\Irr(G),\qquad \chi(1)=d^2.
\]
This character is faithful. Indeed, if $K=\ker\chi$, then
$K\cap N=1$ because $\chi_N=d\theta$. Hence $[K,N]=1$.
The action of $G/N$ on $N/Z(N)\cong B/A$ is faithful, so
$\C_G(N)\le N$, and therefore $K=1$. Since $N$ is nonabelian
and $N/Z(N)\cong C_d^2$, its order is divisible by $d^3$.
Also $|\Sp_2(d)|_d=d$, so $d^4\mid |G|$. It follows that
$d^2\mid\cod(\chi)$, contrary to \eqref{eq:sf}. Thus
$\dl(G)\le d+1$, as required.
\end{proof}

\begin{lemma}\label{lem:third-derived}
Let $r$ be a prime, and let $Y\le\GL_d(\mathbb F_r)$ be solvable and
completely reducible, where $d\le3$. If $Y$ satisfies \eqref{eq:sf},
then $Y^{(3)}$ is an abelian $3$-group.
\end{lemma}

\begin{proof}
As in the preceding proof, it suffices to consider irreducible $Y$.
In degree at most $2$, Lemma~\ref{lem:linear-derived} gives
$Y^{(3)}=1$. In degree $3$, the same holds if $Y$ is imprimitive,
since the kernel on the three one-dimensional blocks is abelian and
the permutation image is contained in $S_3$.

Suppose that $Y$ is primitive of degree $3$, and choose a maximal
solvable overgroup $H\le\GL_3(\mathbb F_r)$. Let $A\le B\le C$
and $a,e$, with $3=ae$, be as in the preceding proof and
\cite[Lemma~2.2]{Seress}. The case $e=1$ is metabelian.
Otherwise $A$ consists of
scalars, $B/A\cong C_3^2$, and $H/B\le\SL_2(3)$.
Consequently $Y^{(3)}\le B\cap Y$. Every $3'$-element of $B$ lies
in $A$. If such an element $cI$ also lies in $Y'$, then
\[
 1=\det(cI)=c^3,
\]
so $c=1$. Thus $B\cap Y'$ contains no nontrivial $3'$-element,
and $Y^{(3)}$ is a $3$-group. It is abelian by
Lemma~\ref{lem:linear-derived}.
\end{proof}

For a $p$-group $P$, let $P^p=\langle x^p:x\in P\rangle$.
Commutation induces an $\Aut(P)$-equivariant epimorphism
\begin{equation}\label{eq:commutator-map}
 \begin{split}
 \beta_P:\bigwedge^2(P/\PhiG(P))&\longrightarrow
           P'/((P')^p[P',P]),\\
 x\PhiG(P)\wedge y\PhiG(P)&\longmapsto [x,y](P')^p[P',P].
 \end{split}
\end{equation}
This is well defined because $\PhiG(P)=P^pP'$, and commutators
are central of exponent $p$ modulo $(P')^p[P',P]$. If $P$ is
nonabelian, its target is nonzero. Indeed, otherwise the finite
abelian $p$-group $P'/[P',P]$ would equal its subgroup of $p$th
powers, so $P'=[P',P]$; nilpotence would then give $P'=1$.

\begin{lemma}\label{lem:primitive-five}
Let $V$ be a nonzero vector space over a prime field $\mathbb F_q$,
and let $L$ act faithfully and irreducibly on $V$. If $V\rtimes L$
is solvable and satisfies \eqref{eq:sf}, then $\dl(L)\le5$.
\end{lemma}

\begin{proof}
Since $L\cong(V\rtimes L)/V$, Lemma~\ref{lem:inheritance}
shows that $L$ and all its subnormal sections satisfy \eqref{eq:sf}.

Let $E=\F(L)$ and $Q=L/E$. The group $E/\PhiG(L)$ is a faithful
completely reducible $Q$-module by \cite[Theorem~1.12]{ManzWolf}.
Each of its prime components is a quotient of
$\OO_p(L)/\PhiG(\OO_p(L))$. Lemma~\ref{lem:affine-rank}
therefore shows that every irreducible constituent has dimension
at most $3$. The group $Q$ embeds in the product of its images on
these constituents, and every image satisfies \eqref{eq:sf}.
Lemmas~\ref{lem:linear-derived} and \ref{lem:third-derived} give
\[
\dl(Q)\le4,\qquad Q^{(3)}\text{ is an abelian }3\text{-group}.
\]

Let $M=L^{(3)}E$, the inverse image of $Q^{(3)}$ in $L$. Then
\[
 M/E=Q^{(3)},\qquad L^{(4)}\le M'\le E.
\]
It is enough to show that every Sylow subgroup of $L^{(4)}$ is abelian.
Let $E_{3'}$ be the Hall $3'$-subgroup of $E$. The group
$M/E_{3'}$ is a normal $3$-subgroup of $L/E_{3'}$.
Lemma~\ref{lem:inheritance} and Lemma A show that it is metabelian, so
$M''\le E_{3'}$. Hence, for $D=L^{(4)}\cap\OO_3(E)$,
\[
 D'\le L^{(5)}=(L^{(3)})''\le M''\le E_{3'}.
\]
As $D'$ is a $3$-group, it is trivial.

Now let $p\ne3$ and $P=\OO_p(E)$. We may assume $P\ne1$.
Let $S$ be a Sylow $3$-subgroup of $M$, and define
\[
 R=P'[P,M]=P'[P,S].
\]
The equality follows from $M=ES$, since the other Sylow subgroups
of $E$ centralize $P$. Both $P$ and $M$ are normal in $L$, so
$R\lhd L$. Moreover,
$M'=E'[E,S]S'$, and $E'[E,S]=[E,M]\le E$. Thus
$M'/[E,M]$ is a $3$-group, whence
\begin{equation}\label{eq:R-PS}
 L^{(4)}\cap P\le M'\cap P\le [E,M]\cap P=R.
\end{equation}
Suppose, for a contradiction, that $R$ is nonabelian. If $S$ acts trivially on
$R/\PhiG(R)$, then it centralizes $R$, since the kernel of
$\Aut(R)\to\Aut(R/\PhiG(R))$ is a $p$-group. The coprime action
identity $[P,S,S]=[P,S]$ would give $[P,S]=1$, so $R=P'$.
This is abelian by Lemma A and Lemma~\ref{lem:inheritance}, a
contradiction. Thus $S$ acts nontrivially on $R/\PhiG(R)$.
Since $p\ne3$, this action is completely reducible. It follows
that $S$ acts nontrivially on some irreducible $L$-composition
factor $U$ of $R/\PhiG(R)$.

The normal $p$-subgroup $P$ acts trivially on $U$, and the other
Sylow subgroups of $E$ centralize $R$. Thus $E$ acts trivially
on $U$. Let $Y$ be the image of $L$ on $U$. The image of
$M=L^{(3)}E$ is $Y^{(3)}$, so $Y^{(3)}\ne1$.
Applying Lemma~\ref{lem:affine-rank} to the normal $p$-subgroup
$R\lhd L$ gives $d(R)\le4$, and every irreducible
$L$-composition factor of $R/\PhiG(R)$ has dimension at most $3$.
Thus $\dim U\le3$. Lemma~\ref{lem:linear-derived}
excludes $\dim U\le2$, since $Y^{(3)}\ne1$. Hence
\[
\dim U=3,\qquad \dl(Y)=4.
\]

Suppose first that $R$ has an abelian normal subgroup $A$ of index
$p$. Lemma~\ref{lem:affine-rank} gives $d(R)\le3$, so
$U=R/\PhiG(R)$ is irreducible. Since $\PhiG(R)\le A$, the
subspace $A/\PhiG(R)$ has dimension $2$. As $A$ is abelian,
\[
 0\ne\bigwedge^2(A/\PhiG(R))\le\ker\beta_R,
\]
where $\beta_R$ is the map in \eqref{eq:commutator-map}.
But $\bigwedge^2U\cong U^*\otimes\bigwedge^3U$ is irreducible,
so the nonzero invariant kernel forces $\beta_R=0$, contrary to
the nonabelianness of $R$.

Lemma A now gives $|R:Z(R)|=p^3$. If $d(R)=3$, then
$U=R/\PhiG(R)$ is irreducible. The image of $Z(R)$ in $U$ is an
invariant subspace. If it were nonzero, irreducibility would give
$Z(R)\PhiG(R)=R$. The Burnside basis theorem would then imply
that $R$ is generated by central elements, a contradiction.
Hence $Z(R)\le\PhiG(R)$, and comparison of indices gives
$Z(R)=\PhiG(R)$.

If $d(R)=4$, let $C=Z(R)\PhiG(R)/\PhiG(R)$. This is a nonzero
proper invariant subspace of $R/\PhiG(R)$, as in the proof of
Lemma~\ref{lem:affine-rank}. It cannot have dimension $2$, since
then all composition factors would have dimension at most $2$.
It cannot have dimension $3$, since three central representatives
together with one further element would generate $R$, making it
abelian. Hence $\dim C=1$, and $(R/\PhiG(R))/C$ is isomorphic
to the $3$-dimensional irreducible factor $U$. Since
$|R:Z(R)|=p^3$, we obtain $\PhiG(R)\le Z(R)$ and
$R/Z(R)\cong U$.
Thus in either case
\[
 \PhiG(R)\le Z(R),\qquad R/Z(R)\cong U.
\]
In particular, $R$ has class $2$, and $R/Z(R)$ is an elementary
abelian group of order $p^3$. For $x,y\in R$,
\[
 [x,y]^p=[x^p,y]=1,
\]
so $R'$ is elementary abelian. Identifying $U$ with $R/Z(R)$,
commutation induces a nonzero equivariant epimorphism
\[
 \bigwedge^2U\longrightarrow R',\qquad
 (xZ(R))\wedge(yZ(R))\longmapsto[x,y].
\]
Its source is irreducible of dimension $3$, whence $R'\cong C_p^3$.

The action of $R$ on $V^*$ is faithful, and $p\ne q$ by
Lemma~\ref{lem:affine-rank}. Let $k$ be a splitting field for $R'$
of characteristic $q$. Choose a nontrivial homomorphism
$\lambda:R'\to k^\times$ whose eigenspace on
$k\otimes_{\mathbb F_q}V^*$ is nonzero, and let $K=\ker\lambda$.
Then $|K|=p^2$. Fixed spaces commute with extension of scalars, so
$\C_{V^*}(K)\supsetneq\C_{V^*}(R')$. Choose
$0\ne v\in\C_{V^*}(K)$ not fixed by $R'$. Then
\[
 |\C_R(v)|\ge p^2,\qquad R'\nleq\C_R(v).
\]
Every subgroup of index $p$ in $R$ contains $R'$. Thus
$|R:\C_R(v)|\ge p^2$, contrary to \eqref{eq:orbit} applied to
$R$. This proves that $R$ is abelian.

By \eqref{eq:R-PS}, every Sylow subgroup of $L^{(4)}$ is abelian.
Since $L^{(4)}\le E$ is nilpotent, $L^{(4)}$ is abelian, and
$\dl(L)\le5$.
\end{proof}

\begin{proof}[Proof of Theorem C]
Let $F=\F(G)$ and $\overline G=G/F$. By
\cite[Theorem~1.12]{ManzWolf}, the group $F/\PhiG(G)$ is a
faithful completely reducible $\overline G$-module, possibly of
mixed characteristic.
Write
\[
 F/\PhiG(G)=V_1\oplus\cdots\oplus V_t
\]
as a direct sum of irreducible $\overline G$-modules, and set
\[
 C_i=\C_{\overline G}(V_i),
 \qquad L_i=\overline G/C_i.
\]
Since $\bigcap_iC_i=1$, the group $\overline G$ is a subdirect subgroup
of $\prod_iL_i$.

In $G/\PhiG(G)$, choose a complement $K$ to $F/\PhiG(G)$.
For each $i$, let $\pi_i:F/\PhiG(G)\to V_i$ be the
$K$-equivariant projection onto $V_i$. The map
\[
 \begin{aligned}
 (F/\PhiG(G))\rtimes K&\longrightarrow
 V_i\rtimes\bigl(K/\C_K(V_i)\bigr),\\
 (w,k)&\longmapsto\bigl(\pi_i(w),k\C_K(V_i)\bigr)
 \end{aligned}
\]
is a surjective homomorphism. Since $K/\C_K(V_i)\cong L_i$,
the group $V_i\rtimes L_i$ is a quotient of $G$.
Hence it satisfies \eqref{eq:sf}, and
Lemma~\ref{lem:primitive-five} gives $\dl(L_i)\le5$. Therefore
\[
 \dl(G/F)=\dl(\overline G)\le5,\qquad G^{(5)}\le F.
\]

Suppose that $G^{(6)}\ne1$. Then $D=G^{(5)}$ is a nonabelian subgroup of
the nilpotent group $\F(G)$. Hence some Sylow $p$-subgroup $P$ of $D$ is
nonabelian. Let $K$ be the Hall $p'$-subgroup of $D$. Since $K$ is
characteristic in $D$, it is normal in $G$. Replacing $G$ by $G/K$ and
retaining the same notation, we may assume that $G^{(5)}=P$ is a
nonabelian $p$-group. Set
\[
 H=G^{(4)}.
\]
Then $H'=P$.

Set $N=\OO_p(H)$. Since $P=H'\le N$, the quotient $H/N$ is abelian. It
has order prime to $p$. For otherwise the inverse image in $H$ of the
Sylow $p$-subgroup of $H/N$ would be a normal $p$-subgroup properly
containing $N$. By Schur--Zassenhaus,
\[
 H=N\rtimes A
\]
for an abelian $p'$-group $A$. Therefore
\begin{equation}\label{eq:last-derived}
 P=H'=N'[N,A].
\end{equation}
The subgroup $N$ is characteristic in $H$ and hence normal in $G$. By
Lemma~\ref{lem:inheritance}, it satisfies \eqref{eq:sf}, and
Lemma A shows that it is metabelian. It is nonabelian, since otherwise
\eqref{eq:last-derived} would make $P$ abelian.

Since $N$ is nonabelian, Lemma A leaves two possibilities. We first
consider the case in which $N$ has an abelian normal subgroup $B$ of
index $p$.
Let
\[
 V=N/\PhiG(N).
\]
Since $B$ is a maximal subgroup, $\PhiG(N)\le B$. Choose
$0\ne f\in V^*$ with
\[
 B/\PhiG(N)=\ker f.
\]
Let $\beta_N$ be the nonzero $G$-equivariant commutator map
in \eqref{eq:commutator-map}, with source $\bigwedge^2V$.

For $g\in G$, let $f^g$ denote the image of $f$ under the dual
action of $G$ on $V^*$. Since $B^g$ is abelian,
\[
 \bigwedge^2\ker(f^g)\le\ker\beta_N.
\]
Set
\[
 W=\langle f^g:g\in G\rangle\le V^*.
\]
We have $\dim W\le2$. Otherwise choose three linearly independent
members $f_1,f_2,f_3$ of the orbit of $f$. Then
\[
 \bigwedge^2V=
 \bigwedge^2\ker f_1+
 \bigwedge^2\ker f_2+
 \bigwedge^2\ker f_3.
\]
Indeed, after extending $f_1,f_2,f_3$ to a basis of $V^*$, every
exterior basis vector lies in at least one of the three summands. Thus
$\beta_N=0$, a contradiction. Hence
\begin{equation}\label{eq:dual-two}
 \dim W\le2.
\end{equation}

The group $A$ acts nontrivially on $W$. Otherwise $A$ fixes $f$, so it
normalizes $B$ and acts trivially on $N/B$. Hence $[N,A]\le B$. Also
$N'\le B$, because $N/B$ is cyclic. Equation~\eqref{eq:last-derived}
would then give $P\le B$, contrary to the nonabelian choice of $P$.

Thus $A$ acts nontrivially on some irreducible $G$-composition factor
$U$ of $W$. If $A$ induced the identity on every such factor, then its
image on $W$ would be a $p$-group. The same image is a quotient of the
$p'$-group $A$, and hence would be trivial, a contradiction. Let $Y$ be the image of $G$ on $U$. By
\eqref{eq:dual-two}, $\dim U\le2$, while $A\le G^{(4)}$ has nontrivial
image on $U$. Thus $Y^{(4)}\ne1$. On the other hand, $Y$ is a quotient
of $G$ and is a completely reducible linear group of degree at most $2$,
so Lemma~\ref{lem:linear-derived} gives $\dl(Y)\le3$, a contradiction.

We are left with
\[
 |N:Z(N)|=p^3.
\]
The group $A$ acts nontrivially on $Q=N/Z(N)$. Otherwise
$[N,A]\le Z(N)$, and \eqref{eq:last-derived} would express $P$ as the
product of the abelian group $N'$ with the central subgroup $[N,A]$;
then $P$ would be abelian.

Set
\[
 \overline Q=Q/\PhiG(Q).
\]
The action of $A$ on $\overline Q$ is still nontrivial. The kernel of
\[
 \Aut(Q)\longrightarrow\Aut(Q/\PhiG(Q))
\]
is a $p$-group, whereas $A$ is a $p'$-group. Moreover,
\[
 \dim_{\mathbb F_p}\overline Q\le3.
\]
Both $Z(N)$ and the inverse image in $N$ of $\PhiG(Q)$ are characteristic
in $N$, so $\overline Q$ is a $G$-module. The nontrivial action of $A$
survives on some irreducible $G$-composition factor $U$ of $\overline Q$:
otherwise the image of $A$ on $\overline Q$ would be a $p$-group and also
a quotient of the $p'$-group $A$, hence trivial. If $Y$ is the image of
$G$ on $U$, then
\[
 \dim U\le3,
 \qquad Y^{(4)}\ne1.
\]
This again contradicts Lemma~\ref{lem:linear-derived}, which gives
$\dl(Y)\le4$.

Both nonabelian alternatives are impossible. Therefore every Sylow
subgroup of $G^{(5)}$ is abelian. Since $G^{(5)}\le\F(G)$ is nilpotent,
$G^{(5)}$ is abelian, and hence $G^{(6)}=1$.
\end{proof}

\section{Sharpness and examples}\label{sec:examples}

For every $\chi\in\Irr(G)$,
\[
 \chi(1)\cod(\chi)=|G:\ker\chi|.
\]
Consequently, if $p^4\nmid |G|$, then
$p^2\nmid\gcd(\chi(1),\cod(\chi))$ for every $\chi\in\Irr(G)$.

\begin{proposition}\label{prop:height-sharp}
There exists a solvable group $G$ satisfying \eqref{eq:sf} with
$h(G)=\dl(G)=5$.
\end{proposition}

\begin{proof}
Let $V=\mathbb F_{11}^4$, and choose $\omega\in\mathbb F_{11}^{\times}$
of order $5$. Define
\[
 D=\left\{\operatorname{diag}(\omega^{a_1},\ldots,\omega^{a_4}):
 a_i\in\mathbb F_5,\quad a_1+\cdots+a_4=0\right\}.
\]
Then $D\cong C_5^3$. The permutation matrices form a subgroup
$S\cong S_4$ of $\GL(V)$ that normalizes $D$ and acts faithfully on it.
Let $H=DS=D\rtimes S$ and $G=V\rtimes H$, with the natural matrix action
on $V$.

Since $\C_H(D)=D$ and $5\nmid|H/D|$, we have $\F(H)=D$.
Similarly, $\C_G(V)=V$ and $11\nmid|G/V|$, so $\F(G)=V$.
Indeed, in each case a Sylow subgroup of the Fitting subgroup belonging
to another prime would centralize the displayed normal Sylow subgroup.
Since $\F(S_4)\cong C_2^2$ and $S_4/\F(S_4)\cong S_3$, it follows that
\[
 h(G)=1+h(H)=2+h(S_4)=5.
\]
The groups $V$ and $D$ are abelian and $\dl(S_4)=3$, so
$\dl(G)\le5$. The inequality $h(G)\le\dl(G)$ gives $\dl(G)=5$.

We have $|G|=11^4\cdot5^3\cdot2^3\cdot3$.
Since $V$ is abelian and normal, \cite[Theorem~6.15]{Isaacs/book}
shows that every irreducible character degree of $G$ divides $|G:V|$; in particular, no such degree is divisible by $11$.
For each remaining prime $p\mid |G|$, we have $p^4\nmid|G|$.
Hence \eqref{eq:sf} holds.
\end{proof}

The next construction attains the derived-length bound.

\begin{proposition}\label{prop:derived-sharp}
There exists a solvable group $G$ satisfying \eqref{eq:sf} with
$\dl(G)=6$ and $h(G)=3$.
\end{proposition}

\begin{proof}
Let $E$ be the extraspecial group of order $3^3$ and exponent $3$.
To specify the action used below, realize $E$ on
$\mathbb F_3^2\oplus\mathbb F_3$ with multiplication
\[
 (u,a)(v,b)=\left(u+v,a+b+\tfrac12\langle u,v\rangle\right),
\]
where $\langle\ ,\ \rangle$ is the standard nondegenerate
alternating form on $\mathbb F_3^2$. The group $\SL_2(3)$ acts on
$E$ by $g(u,a)=(gu,a)$. Let $Q_8$ be its quaternion subgroup of
order $8$, and let $H=E\rtimes Q_8$. Thus $Q_8$ fixes $Z(E)$
and acts faithfully on $E/Z(E)$.

Choose a nonlinear $\theta\in\Irr(E)$. It has degree $3$, is
faithful, and is determined by its nontrivial central character.
It is therefore $H$-invariant. Since $|H:E|$ is coprime to $|E|$,
it extends to $H$ by \cite[Corollary~6.28]{Isaacs/book}.
Such an extension is faithful.
Indeed, its kernel intersects $E$ trivially and therefore
centralizes $E$. Its image in $Q_8$ acts trivially on $E/Z(E)$,
so that image is trivial; faithfulness on $E$ gives the result.

Let $p=433$. Since $|H|=216$ divides $p-1$, the field
$\mathbb F_p$ is a splitting field for $H$, of characteristic
coprime to $|H|$. The extension above therefore gives a faithful
$3$-dimensional $\mathbb F_pH$-module $W$. The group $Z(E)$
acts on $W$ as a nontrivial scalar group of order $3$.
Define a group $P$ on $W\oplus\bigwedge^2W$ by
\[
 (u,a)(v,b)=\left(u+v,a+b+\tfrac12u\wedge v\right).
\]
Then $P$ has class $2$ and exponent $p$, with
\[
 P'=Z(P)=\PhiG(P)=\bigwedge^2W,\qquad P/Z(P)\cong W.
\]
The action of $H$ on $W$ induces an action on $P$. Let
$G=P\rtimes H$.

For every $\eta\in\Irr(P)$, the inequality
$\eta(1)^2\le|P:Z(P)|=p^3$ gives $\eta(1)\in\{1,p\}$.
If $\chi\in\Irr(G\mid\eta)$, Clifford theory gives, for some
positive integer $e$,
\[
 \chi(1)=e|H:I_H(\eta)|\eta(1),
 \qquad e^2\le|I_H(\eta)|\le216<p.
\]
Hence $v_p(\chi(1))\le1$. The only other primes dividing
$|G|$ are $2$ and $3$, each with exponent $3$. The observation
at the start of this section proves \eqref{eq:sf}.

The nontrivial element of $Q_8'=Z(Q_8)$ acts as $-I$ on
$E/Z(E)$. The Burnside basis theorem gives
$[E,Q_8']=E$, and hence also $[E,Q_8]=E$. Therefore
\[
 H'=E\rtimes C_2,\qquad H''=E,\qquad H'''=Z(E),
 \qquad H^{(4)}=1.
\]
The group $Z(E)$ acts on $W$ as a nontrivial scalar of order $3$,
and on $\bigwedge^2W$ as its square. Thus $[P,Z(E)]=P$, and
the derived series of $G$ is
\[
 \begin{aligned}
 G'&=P\rtimes H',& G''&=P\rtimes H'',&
 G'''&=P\rtimes H''',\\
 G^{(4)}&=P,& G^{(5)}&=P',& G^{(6)}&=1.
 \end{aligned}
\]
In particular, $\dl(G)=6$.

The action of $H$ on $P/\PhiG(P)\cong W$ is faithful, so
$\C_G(P)\le P$. Every Sylow subgroup of $\F(G)$ belonging to
a prime other than $p$ centralizes $P$ and is therefore trivial.
Since $P$ is a normal Sylow $p$-subgroup, $\F(G)=P$.
Similarly, a normal $2$-subgroup of $H$ centralizes $E$ and
therefore acts trivially on $E/Z(E)$, so it is trivial.
Thus $\F(H)=E$, and $h(H)=2$ because $H/E\cong Q_8$.
Consequently $h(G)=3$.
\end{proof}

The rank of $G/\F(G)$ need not be bounded under \eqref{eq:sf},
even when $G$ is solvable.

\begin{proposition}\label{prop:unbounded-rank}
For each positive integer $n$, the group $G=S_3^n$ satisfies
\eqref{eq:sf}, while $G/\F(G)\cong C_2^n$ has rank $n$.
\end{proposition}

\begin{proof}
Write $G=G_1\times\cdots\times G_n$, where $G_i\cong S_3$.
The irreducible characters of $S_3$ are the principal character $1$,
the sign character $\varepsilon$, and the character $\varphi$ of degree $2$.
Every $\chi\in\Irr(G)$ is an external tensor product of
irreducible characters of the factors. Suppose that $a$ factors contribute $\varphi$ and $b$
factors contribute $\varepsilon$. Then $\chi(1)=2^a$.

The degree-$2$ representation of $S_3$ is faithful, and no nonidentity
 element acts as a scalar. Thus every coordinate belonging to one of the
 $a$ degree-$2$ factors must be the identity for an element of $G$ to lie
 in $\ker\chi$. Indeed, a tensor product of invertible matrices can be the identity
 only if every factor is scalar. These coordinates therefore contribute
 a factor $6^a$ to $|G:\ker\chi|$.

 For the $b$ sign factors, consider the homomorphism
 \[
 \lambda:S_3^b\longrightarrow\{1,-1\},\qquad
 (g_1,\ldots,g_b)\longmapsto
 \varepsilon(g_1)\cdots\varepsilon(g_b).
 \]
 If $b>0$, this map is onto, so its kernel has index $2$, independent
 of $b$. Explicitly, $\ker\lambda$ consists of the tuples having an
 even number of odd permutations. The factors carrying the principal
 character impose no restriction. Therefore the kernel of $\chi$ has
 index $6^a$ if $b=0$ and $2\cdot6^a$ if $b>0$.
 Consequently,
\[
 |G:\ker\chi|=\begin{cases}
  6^a,&b=0,\\
  2\cdot 6^a,&b>0.
 \end{cases}
\]
It follows that
\[
 \cod(\chi)=\begin{cases}
  3^a,&b=0,\\
  2\cdot 3^a,&b>0.
 \end{cases}
\]
In particular $\gcd(\chi(1),\cod(\chi))$ is either $1$ or $2$.
Finally, $\F(S_3)=C_3$ and therefore
$\F(G)=C_3^n$ and $G/\F(G)\cong C_2^n$.
\end{proof}

\section*{Acknowledgements}
Y. Yang was supported in part by Simons Foundation Award \#918096.

This paper extends K. Aldahleh's honors thesis at the University of
Rochester \cite{Aldahleh}, written under the supervision of Amanda
Tucker. The authors thank her for her guidance and helpful suggestions.

\section*{Disclosure Statement}
The authors report no conflict of interest.

\section*{Data Availability}
No data were used in this study.

\end{document}